\documentclass[twoside,11pt,leqno]{article}
\usepackage{amsfonts}

   \def\H{{\cal H}}

\def\B{B({\cal H})}

\newtheorem{df}{Definition}[section]
\newtheorem{thm}[df]{Theorem} \newtheorem{pro}[df]{Proposition}
\newtheorem{cor}[df]{Corollary} 
 \newtheorem{lem}[df] {Lemma}
\def\sfstp{{\hskip-1em}{\bf.}{\hskip1em}}

\def\subject#1{\renewcommand{\thefootnote}{}\footnote
{AMS(MOS) subject classification (2010). Primary: {#1}}}

\def\keywords#1{\renewcommand{\thefootnote}{}\footnote
{Keywords: {#1}}}

\def\enddemo{\qed \endtrivlist} \expandafter\let\csname
enddemo*\endcsname=\enddemo

\def\qedsymbol{\ifmmode\bgroup\else$\bgroup\aftergroup$\fi
\vcenter{\hrule\hbox{\vrule
height.5em\kern.5em\vrule}\hrule}\egroup}
\def\qed{\ifmmode\else\unskip\nobreak\fi\quad\qedsymbol}

\title{\bf Operator roots of polynomials: iso-symmetric operators}
\author{\normalsize B.P.~Duggal, I.H.~Kim}
\date{}

\begin{document}

\maketitle \thispagestyle{empty} \vskip-16pt

\subject{47A05, 47A55; Secondary47A11, 47B47.} \keywords{ Hilbert space,  Left/right multiplication operator,  $m$-left invertible, $m$-isometric and $m$-selfadjoint operators, product of operators, perturbation by nilpotents, commuting operators.  }
\footnote{The second named author was supported by Basic Science Research Program through the National Research Foundation of Korea (NRF)
 funded by the Ministry of Education (NRF-2019R1F1A1057574)}
\begin{abstract} Given Hilbert space operators $A_i, B_i$, $i=1,2$, and $X$ such that $A_1$ commutes with $A_2$ and $B_!$ commutes with $B_2$, and integers $m, n\geq 1$, we say that the pairs of operators $(B_1,A_1)$ and $(B_2,A_2)$ are left-$(X, (m,n))$-symmetric, denoted $((B_1,A_1),(B_2,A_2))\in {\rm left}-(X,(m,n))-{\rm symmetric}$ if
$$ \sum_{j=0}^m\sum_{k=0}^n (-1)^{j+k}\left(\begin{array}{clcr}m\\j\end{array}\right) \left(\begin{array}{clcr}n\\k\end{array}\right) B_1^{m-j}B_2^{n-k} X A_2^{n-k}A_1^{j}=0.$$An important class of left-$(X,(m,n))-$symmetric operators is obtained uponchoosing $B_1=B_2=A^*_1=A^*_2=A^*$ and $X=I$: such operators have been called $(m,n)-$isosymmetric, and a study of the spectral picture and maximal invariant subspaces of $(m,n)-$isosymmetric operators has been carried out by Stankus \cite{St}. The current work considers  stability under perturbations by commuting nilpotents, and  products of commuting, left-$(X, (m,n))-$symmetric operators. It is seen that  $(X, (m,n))-$isosymmetric Drazin invertible operators $A$ have a particularly interesting structure.
\end{abstract}


\section {\sfstp Introduction} Let $\B$ denote the algebra of operators, i.e. bounded linear transformations, on an infinite dimensional complex Hilbert space $\H$ into itself. Let $\mathbb C$ denote the complex plane, ${\mathbb C}^d$ the product of $d$ copies of $\mathbb C$ for some integer $d\geq 1$, $\overline{z}$ the conjugate of $z\in\mathbb C$ and ${\bf z}=(z_1,z_2,...,z_d)\in{\mathbb C}^d$. For a given polynomial $P$ in ${\mathbb C}^d$ and a $d$-tuple ${\bf A}$ of commuting operators in $\B$, ${\bf A}$ is a hereditary root of $P$ if $P({\bf A})=0$. Hereditary roots have attracted the attention of a number of researchers in the recent past. Two particular operator classes of hereditary roots which have been studied extensively are those of $m$-symmetric (also called $m$-selfadjoint in the literature) and $m$-isometric operators, where $A\in\B$ is $m$-symmetric (for some integer $m\geq 1$) if
$$\sum_{j=0}^m(-1)^j\left(\begin{array}{clcr}m\\j\end{array}\right){A^*}^{(m-j)}A^j=0$$
and $A\in\B$ is $m$-isometric if
$$\sum_{j=0}^m(-1)^j\left(\begin{array}{clcr}m\\j\end{array}\right){A^*}^{(m-j)}A^{m-j}=0.$$
It is clear that $A$ is $m$-symmetric if and only if it arises as a solution of $P(z)=({\overline{z}}-z)^m=0$, and $A$ is $m$-isometric if and only if it arises as a solution of $P(z)=({\overline{z}}z-1)^m=0$. The class of $m$-symmetric operators was introduced by Helton \cite{Hel} ({\it albeit} not as operator solutions of the polynomial equation $({\overline{z}}-z)^m=0$), who showed that an operator $A$ is $2$-Jordan (i.e., $ A=T+N$ for some self-adjoint $T$ and a $2$-nilpotent $N$ commuting with $T$) if and only if $A$ and $A^*$ are $3$-symmetric. McCullough and Rodman \cite{McR} in their consideration of algebraic and spectral properties of $m$-symmetric operators proved that if an $A\in\B$ is self-adjoint and an $N\in\B$ is an $n$-nilpotent which commutes with $A$, then $A+N$ is $(2n-1)$-symmetric. These operators have since been considered by many other authors, amongst them Stankus \cite {St} and Trieu Le \cite{TL}

\

The class of $m$-symmetric operators was introduced by Agler and studied in a series of papers by Agler and Stankus \cite{{AS1}, {AS2}, {AS3}}; properties of $m$-isometric operators, amongst them the spectral picture, strict $m$-isometries, perturbation by commuting nilpotents and the product of $m$-isometries, have since been studied by a large number of authors, amongst them Sid Ahmed \cite {OAM}, Bayart \cite{FB}, Bermudez {\it et al} \cite{{BMN}, {BMN1}, {BMMN}}, Botelho and Jamison \cite{Fb}, Duggal \cite{{BD1}, {BD2}, {BD3}}, Gu \cite{{G}, {G1}} and Gu and Stankus \cite{GS}.

\

$A\in\B$ is an $(m,n)$-isosymmetry for some integers $m,n\geq 1$ if
\begin{eqnarray*} & &  \sum_{j=0}^m(-1)^j\left(\begin{array}{clcr}m\\j\end{array}\right){A^*}^{(m-j)}\left(\sum_{k=0}^n(-1)^k\left(\begin{array}{clcr}n\\k\end{array}\right){A^*}^{(m-k)}A^k \right)A^{m-j}\\
&=& \sum_{k=0}^n(-1)^k\left(\begin{array}{clcr}n\\k\end{array}\right){A^*}^{(m-k)}\left(\sum_{j=0}^m(-1)^j\left(\begin{array}{clcr}m\\j\end{array}\right){A^*}^{(m-j)}A^{m-j}\right) A^k\\
&=& 0; \end{eqnarray*}
$(m,n)$-isosymmetric operators arise as the hereditary roots of the polynomial $({\overline{z}}z-1)^m({\overline{z}}-z)^n=0$, and a study, amongst other properties, of the spectrum, resolvent inequalities and maximal invariant subspaces of these operators has been carried out by Stankus \cite{St}. In this paper we study  a generalisation of $(m,n)$-isosymmetric operators, but from the point of view of elementary operators. The problem that we consider is that of the permanence of this generalised isosymmetric property under commuting products and perturbation by commuting nilpotents.

\

For $A,B\in\B$, let $L_A$ and $R_B \in B(\B)$ denote respectively the operators
$$
L_A(X)=AX \ {\rm and}\ R_B(X)=XB
$$
of left multiplication by $A$ and right multiplication by $B$. We say that the operator $A$ is left $(X, m)$-invertible by $B$, denoted $(B,A)\in$ left-$(X, m)$-invertible, for some operator $X\in\B$  if
$$
\triangle_{B,A}^m(X)=\left(L_BR_A-I\right)^m(X)=\sum_{j=0}^m(-1)^j\left(\begin{array}{clcr}m\\j\end{array}\right)B^{m-j}XA^{m-j}=0,
$$
and that the operator $B$ is an $(X,n)$-symmetry of $A$, denoted $(B,A)\in (X,n)$-symmetry, for some $X\in\B$ if

$$
\delta_{B,A}^n(I)=(L_B-R_A)^n(X)=\sum_{j=0}^n(-1)^j\left(\begin{array}{clcr}n\\j\end{array}\right){B}^{(n-j)}XA^j=0.
$$
It is clear from these definitions that $\triangle_{B,A}^m(I)$ defines the class of left-$m$-invertible operators  $A$ of \cite{DM}, $\triangle^m_{A^*,A}(I)$ defines the class of $m$-isometric operators $A$, 
 $\delta^n_{A^*,A}(I)$ defines the class of $n$-symmetric (equivalently, $n$-selfadjoint) operators, and an
operator $A\in\B$ is $(m,n)$-isosymmetric if and only if
$$\triangle^m_{A^*,A}\left(\delta^n_{A^*,A}(I)\right)= \delta_{A^*,A}^n\left(\triangle_{A^*,A}^m(I)\right)=0.$$
Let $[A, B]=AB-BA$ denote the commutator of $A, B\in\B$. Given operators $A_i, B_i, X\in\B$, $i=1,2$ and  positive integers $m$ and $n$, such that  $$[A_1, A_2]=[B_1, B_2]=0,$$
we say in the following that the pairs of operators $(B_1,A_1)$ and $(B_2,A_2)$ are left-$(X, (m,n))$-symmetric, denoted $$((B_1,A_1),(B_2,A_2))\in {\rm left}-(X,(m,n))-{\rm symmetric}$$ if
\begin{eqnarray*} & & \triangle^m_{B_1,A_1}\left(\delta^n_{B_2,A_2}\right)(X)=(L_{B_1}R_{A_1}-I)^m\left((L_{B_2}-R_{A_2})^n(X)\right)\\
&=& \sum_{j=0}^m\sum_{k=0}^n {(-1)^{j+k}\left(\begin{array}{clcr}m\\j\end{array}\right) \left(\begin{array}{clcr}n\\k\end{array}\right) L^{m-j}_{B_1}L^{n-k}_{B_2}R^k_{A_2}R^{m-j}_{A_1}(X)}\\
&=&   (L_{B_2}-R_{A_2})^n\left((L_{B_1}R_{A_1}-I)^m(X)\right)= \delta^n_{B_2,A_2}\left(\triangle^m_{B_1,A_1}(X)\right)\\
&=&  \sum_{k=0}^n\sum_{j=0}^m {(-1)^{j+k}\left(\begin{array}{clcr}m\\j\end{array}\right) \left(\begin{array}{clcr}n\\k\end{array}\right) L^{n-k}_{B_2}L^{m-j}_{B_1}R^{m-j}_{A_1}R^{k}_{A_2}(X)}\\
&=& 0.
\end{eqnarray*}
Products, and perturbation by commuting nilpotents, of left-$(X,(m,n))$-symmetric pairs of operators behave in a manner very similar to that of $m$-isometric and $n$-symmetric operators. We prove:
\begin{thm}\label{thm1} If $A_i, B_i, S_i, T_i, X\in\B$, $i=1,2$, are such that
\vskip4pt\noindent (i) $[A_1, A_2]=[B_1, B_2]=[A_i, T_i]=[B_i, S_i]=0$,

\vskip4pt\noindent (ii)  $((B_1,A_1),(B_2,A_2))\in$ {\rm left}-$(X,(m_1,n_1))$-{\rm symmetric},

\vskip4pt\noindent (iii)  $((S_1,T_1),(B_2,A_2))\in$ {\rm left}-$(X,(r_1,n_2))$-{\rm symmetric},

\vskip4pt\noindent (iv)  $((B_1,A_1),(S_2,T_2))\in$ {\rm left}-$(X,(m_2,s_1))$-{\rm symmetric}, and

\vskip4pt\noindent (v)  $((S_1,T_1),(S_2,T_2))\in$ {\rm left}-$(X,(r_2,s_2))$-{\rm symmetric},

\vskip4pt\noindent then
$$ ((S_1B_1,T_1A_1), (S_2B_2,T_2A_2))\in {\rm left}-(X,(m+r-1,n+s-1))-{\rm symmetric},$$
where $m={\rm max}(m_1,m_2)$, $n={\rm max}(n_1,n_2)$, $r={\rm max}(r_1,r_2)$ and $s={\rm max}(s_1,s_2)$.
\end{thm}

\begin{thm}\label{thm2} If $A_i, B_i, M_i, N_i, X\in \B$, $i=1,2$, are such that

\vskip4pt\noindent (i) $M_i^{m_i}=N_i^{n_i}=0$, $m_i$ and $n_i$ some positive integers ($i=1,2$),

\vskip4pt\noindent (ii) $[A_1, A_2]=[B_1, B_2]=[M_1, M_2]=[N_1, N_2]=[A_i, M_i]=[B_i, N_i]=0$ ($i=1,2$) and

\vskip4pt\noindent (iii) $((B_1,A_1),(B_2,A_2))\in {\rm left}-(X, (m,n))-{\rm symmetric}$,

\vskip4pt\noindent then
$$ ((B_1+N_1,A_1+M_1), (B_2+N_2,A_2+M_2))\in {\rm left}-(X,(m+m_1+n_1-2,n+m_2+n_2-2))-{\rm symmetric}.$$
\end{thm}
 An operator $A\in \B$ is Drazin invertible, with Drazin inverse $A_d$, if
$$
[A_d,A]=0,\ A_d^2A=A_d,\ A^{p+1}A_d=A^p
$$
for some integer $p\geq 1$. (The least integer $p$ for which this holds is then called the Drazin index of $A$.) No Drazin invertible operator $A\in\B$ can be $m$-isometric (equivalently, left-$m$-invertible by its adjoint) \cite{DK}:
there may however exist operators $X\in\B$ such that  $A$ is  left-$(X,m)$-invertible by $A^*$. We prove:

\begin{thm}\label{thm3} Let $A\in\B$ be a Drazin invertible operator, with Drazin index $p$ and Drazin inverse $A_d$. Let $X\in\B$, and let $m,n$ be some positive integers.

\vskip4pt\noindent (i)  If $A\in (X,(m,n))$-isosymmetric, then $$\triangle^n_{A^*_d,A}\left(\triangle^m_{A^*,A_d}(X)\right)=0=\delta^m_{A^*_d,A}\left(\delta^n_{A^*_d,A_d}(X)\right);$$

\vskip4pt\noindent (ii) if $(A^*_d,A),(A^*,A))\in$ left-$(X,(m,n))$-symmetric, then $\triangle^n_{A^*_d,A}\left(\triangle^m_{A^*,A_d}(X)\right)=0$;

\vskip4pt\noindent (iii) if $((A^*,A),(A^*_d,A))\in$ left-$(X,(m,n))$-symmetric, then $$\delta^n_{A^*,A_d}\left(\delta^m_{A^*_d,A}(X)\right)=0.$$
\end{thm}

We prove Theorem \ref{thm1}, and most of our complementary results, in Section 2, tensor products are considered in Section 3, Section 4 is devoted to the proof of Theorem \ref{thm2}, and we prove Theorem \ref{thm3} in Section 5..
\

\section {\sfstp Complementary results, Proof of  Theorem \ref{thm1}.} We start this section by proving some complementary results. Throughout the following $A_i, B_i$, $i=1,2$, and $X$ will denote operators in $\B$, and $m, n, t$ will denote positive integers.

\begin{lem}\label{lem0}  \cite[Proposition 2.1]{DK}  If $(B,A))\in \{(X,n)-symmetric\}\vee \{ left-(X,m)- invertible\}$;
 \vskip4pt\noindent (i) then $(B,A)\in \{(X,t)-{symmetric}\} \vee ({resp.}) \{{left}-(X,t)-{ invertible}\}$ for all integers $t\geq n,m$;
\vskip4pt\noindent (ii)  and $A, B$ are invertible, then $(B^{-1},A^{-1})\in\{(X,n)-{symmetric}\} \vee ({resp.}) \{{left}-(X,m)-{invertible}\}$.
\end{lem}
The following lemma is an easy consequence of Lemma \ref{lem0} (ii).

\begin{lem}\label{lem1} If $A_i, B_i$ are invertible and $((B_1,A_1),(B_2,A_2))\in$ left-$(X,(m,n))$-symmetric, then $((B_1^{-1},A_1^{-1}),(B_2^{-1},A_2^{-1}))\in$ left-$(X,(m,n))$-symmetric.
\end{lem}
Trivially, $[L_B, R_A]=0$ for all operators $A, B\in\B$. Hence, for every integer $k\geq 1$,
\begin{eqnarray*}
(L_B^k-R_A^k)^m &=& (L_{B^k}-R_{A^k})^m=(L_B-R_A)^m P_{m(k-1)}(L_B,R_A)\\ &=&P_{m(k-1)}(L_B,R_A)(L_B-R_A)^m
\end{eqnarray*}
and
\begin{eqnarray*}
(L_B^kR_A^k-I)^m&=&(L_{B^k}R_{A^k}-I)^m=(L_BR_A-I)^mQ_{m(k-1)}(L_B,R_A)\\ &=&Q_{m(k-1)}(L_B,R_A)(L_BR_A-I)^m
\end{eqnarray*}
for some polynomials $P$ and $Q$ of degree $m(k-1)$. We have:
\begin{lem}\label{lem2}  \cite{DK}  If  $(B,A))\in \{(X,n)-{symmetric}\} \vee \{{left}-(X,m)-{invertible}\}$, then  $(B^k,A^k))\in \{(X,n)-{symmetric}\} \vee {(resp.)} \{{left}-(X,m)-{invertible}\}$ for every positive integer $k$.
\end{lem}

Considering operators $A_i$ and $B_i$, $i=1,2$, such that $[A_1, A_2]=[B_1, B_2]=0$, one has
$$
[L_{B_i}, R_{A_j}]=[L_{B_1},L_{B_2}]=[R_{A_1}, R_{A_2}]=0; \ i,j=1,2.
$$
Hence
$$
\triangle^m_{B_1,A_1}\left(\delta^n_{B_2,A_2}(X)\right)=\delta^n_{B_2,A_2}\left(\triangle^m_{B_1,A_1}(X)\right).
$$
We have:
\begin{lem}\label{lem3} If $[A_1, A_2]=[B_1, B_2]=0$, and either $(B_1,A_1)\in$ left-$(X,m)$-{\rm invertible} or $(B_2,A_2)\in (X,n)$- {\rm symmetric}, then $\triangle^m_{B_1,A_1}\left(\delta^n_{B_2,A_2}(X)\right)=\delta^n_{B_2,A_2}\left(\triangle^m_{B_1,A_1}(X)\right)=0$.
\end{lem}

\begin{lem}\label{lem4} If $[A_1, A_2]=[B_1, B_2]=0$ and $((B_1,A_1),(B_2,A_2))\in$ {\rm left}-$(X,(m,n))$-{\rm symmetric}, then $((B_1,A_1),(B_2,A_2))\in \{$ {\rm left}-$(X,(m_1,n))$-{\rm symmetric}$\} \wedge \{$ {\rm left}-$(X,(m,n_1))$-{\rm symmetric}$\} \wedge
\{${\rm left}-$(X,(m_1,n_1))$-{\rm symmetric}$\}$ for all integers $m_1\geq m$ and $n_1\geq n$.
\end{lem}

\begin{demo} The proof follows since
$$
\triangle^{m_1}_{B_1,A_1}(\delta^n_{B_2,A_2})=\triangle^{m_1-m}_{B_1,A_1}\left(\triangle^m_{A_1,B_1}(\delta_{B_2,A_2}^n)\right),
$$
$$
\triangle^m_{B_1,A_1}\left(\delta^{n_1}_{B_2,A_2}\right)=\triangle^m_{B_1,A_1}\left(\delta^{n_1-n}_{B_2,A_2}\delta^n_{B_2,A_2}\right)=\delta^{n_1-n}_{B_2,A_2}
\left(\triangle^m_{B_1,A_1}\delta^{n}_{B_2,A_2}\right)
$$
and
$$
\triangle^{m_1}_{B_1,A_1}\left(\delta^{n_1}_{B_2,A_2}\right)=\triangle^{m_1-m}_{B_1,A_1}\delta^{n_1-n}_{B_2,A_2}\left(\triangle^{m}_{B_1,A_1}\delta^{n}_{B_2,A_2}\right).
$$
\end{demo}

\begin{pro}\label{pro1} Let $A_i, B_i, S,T,X \in \B,\ i=1,2,$ be such that
$$
[A_1,A_2]=[B_1,B_2]=[A_i,T]=[B_i,S]=0,
$$
and let (a), (b), (c), (d) and (e) be the hypotheses:

(a) $\left((B_1,A_1), (B_2,A_2)\right)\in$ left-$\left(X,(m,n)\right)$-symmetric;

(b) $(B_1,A_1)\in$ left-$\left(X,m)\right)$-invertible;

(c) $(B_2,A_2)\in$ left-$\left(X,n)\right)$-symmetric;

(d) $(S,T)\in$ left-$\left(X,t)\right)$-invertible;

(e) $(S,T)\in$ left-$\left(X,t)\right)$-symmetric.

\vskip4pt\noindent (i) If either of the hypotheses (b) or (a)$\wedge$(e) or (c)$\wedge$(e) is satisfied, then
$$
\left((B_1,A_1), (SB_2,TA_2)\right)\in left-\left(X,(m,n+t-1)\right)-symmetric.
$$

\vskip4pt\noindent (ii) If either of the hypotheses (c) or (a)$\wedge$(d) or (b)$\wedge$(d) is satisfied, then
$$
\left((SB_1,TA_1), (B_2,A_2)\right)\in left-\left(X,(m+t-1,n)\right)-symmetric.
$$
\end{pro}

\begin{demo} The hypotheses $[S,B_2]=[T,A_2]=0$ imply
\begin{eqnarray*} \delta_{SB_2,TA_2}^{n+t-1}&=&\left(L_SL_{B_2}-R_T R_{A_2}\right)^{n+t-1}\\
&=& \left\{L_S\left(L_{B_2}-R_{A_2}\right)+\left(L_S-R_T\right)R_{A_2}\right\}^{n+t-1}\\
&=&  \sum_{j=0}^{n+t-1}\left(\begin{array}{clcr}n+t-1\\j\end{array}\right) L_S^{n+t-1-j}R_{A_2}^j\delta_{B_2,A_2}^{n+t-1-j}\delta_{S,T}^j\\
&=& \sum_{j=0}^{n+t-1}\left(\begin{array}{clcr}n+t-1\\j\end{array}\right) L_S^{n+t-1-j}R_{A_2}^j\delta_{S,T}^j\delta_{B_2,A_2}^{n+t-1-j}
\end{eqnarray*}
and the hypotheses $[S,B_1]=[T,A_1]=0$ imply
\begin{eqnarray*} \triangle_{SB_1,TA_1}^{m+t-1}&=&\left(L_S L_{B_1}R_T R_{A_1}-I\right)^{m+t-1}\\
&=& \left\{L_S R_T\left(L_{B_1}R_{A_1}-I\right)+\left(L_S R_T-I\right)\right\}^{m+t-1}\\
&=& \sum_{k=0}^{m+t-1}\left(\begin{array}{clcr}m+t-1\\k\end{array}\right) \left(L_S R_T\right)^{m+t-1-k}\triangle_{B_1,A_1}^{m+t-1-k}\triangle_{S,T}^k\\
&=&\sum_{k=0}^{m+t-1}\left(\begin{array}{clcr}m+t-1\\k\end{array}\right) \left(L_S R_T\right)^{m+t-1-k}\triangle_{S,T}^k\triangle_{B_1,A_1}^{m+t-1-k}.
\end{eqnarray*}
This, in view of the commutativity hypotheses on $A_i,B_i,S$ and $T \ (i=1,2)$, implies
\begin{eqnarray*}
 & & \triangle_{B_1,A_1}^{m}\left(\delta_{SB_2,TA_2}^{n+t-1}(X)\right)\\
  &=&\sum_{j=0}^{n+t-1}\left(\begin{array}{clcr}n+t-1\\j\end{array}\right)L_S^{n+t-1-j} R_{A_2}^j\delta_{S,T}^j\delta_{B_2,A_2}^{n+t-1-j}\left(\triangle_{B_1,A_1}^{m}(X)\right)\\
&=&\sum_{j=0}^{t-1}\left(\begin{array}{clcr}n+t-1\\j\end{array}\right)L_S^{n+t-1-j} R_{A_2}^j\delta_{S,T}^j\left(\triangle_{B_1,A_1}^{m}\delta_{B_2,A_2}^{n+t-1-j}(X)\right)\\
& & +\sum_{j=t}^{n+t-1}\left(\begin{array}{clcr}n+t-1\\j\end{array}\right)L_S^{n+t-1-j} R_{A_2}^j\triangle_{B_1,A_1}^{m}\delta_{B_2,A_2}^{n+t-1-j}\left(\delta_{S,T}^j(X)\right)\\
&=&\sum_{j=0}^{t-1}\left(\begin{array}{clcr}n+t-1\\j\end{array}\right)L_S^{n+t-1-j} R_{A_2}^j\delta_{S,T}^j\triangle_{B_1,A_1}^{m}\left(\delta_{B_2,A_2}^{n+t-1-j}(X)\right)\\
& &+\sum_{j=t}^{n+t-1}\left(\begin{array}{clcr}n+t-1\\j\end{array}\right)L_S^{n+t-1-j} R_{A_2}^j\triangle_{B_1,A_1}^{m}\delta_{B_2,A_2}^{n+t-1-j}\left(\delta_{S,T}^j(X)\right)
\end{eqnarray*}
and
\begin{eqnarray*}
& & \triangle_{SB_1,TA_1}^{m+t-1}\left(\delta_{B_2,A_2}^{n}(X)\right)\\
&=&\sum_{k=0}^{m+t-1}\left(\begin{array}{clcr}m+t-1\\k\end{array}\right)\left(L_S R_T\right)^{m+t-1-k} \triangle_{S,T}^k\triangle_{B_1,A_1}^{m+t-1-k}\left(\delta_{B_2,A_2}^n(X)\right)\\
&=&\sum_{k=0}^{t-1}\left(\begin{array}{clcr}m+t-1\\k\end{array}\right)\left(L_S R_T\right)^{m+t-1-k} \triangle_{S,T}^k\left(\triangle_{B_1,A_1}^{m+t-1-k}\delta_{B_2,A_2}^n(X)\right)\\
& & +\sum_{k=t}^{m+t-1}\left(\begin{array}{clcr}m+t-1\\k\end{array}\right)\left(L_S R_T\right)^{m+t-1-k} \triangle_{B_1,A_1}^{m+t-1-k}\left(\triangle_{B_1,A_1}^k\delta_{B_2,A_2}^n(X)\right)\\
&=&\sum_{k=0}^{t-1}\left(\begin{array}{clcr}m+t-1\\k\end{array}\right)\left(L_S R_T\right)^{m+t-1-k} \triangle_{S,T}^k\left(\delta_{B_2,A_2}^n\triangle_{B_1,A_1}^{m+t-1-k}(X)\right)\\
& & +\sum_{k=t}^{m+t-1}\left(\begin{array}{clcr}m+t-1\\k\end{array}\right)\left(L_S R_T\right)^{m+t-1-k} \triangle_{B_1,A_1}^{m+t-1-k}\delta_{B_2,A_2}^n\left(\triangle_{S,T}^k(X)\right).
\end{eqnarray*}
Recall now that
$$
\triangle_{B_1,A_1}^{m}\delta_{B_2,A_2}^n(X)=0 \Longrightarrow \triangle_{B_1,A_1}^{m}\delta_{B_2,A_2}^{n_1}(X)=0 \ {\rm for \ all}\ n_1 \geq n
$$
and
$$
\delta_{S,T}^t(X)=0 \Longrightarrow \delta_{S,T}^{t_1}(X)=0 \ {\rm for \ all}\ t_1 \geq t
$$
(see Lemmas \ref{lem4} and \ref{lem0}). Hence, if either of the hypotheses (b), or, $(a)\wedge (e)$, or, $(c)\wedge (e)$ is satisfied, then
$$
\left((B_1,A_1), (SB_2,TA_2)\right)\in left-\left(X,(m,n+t-1)\right)-symmetric.
$$
Again, since
$$
\triangle_{B_1,A_1}^{m}\left(\delta_{B_2,A_2}^{n}(X)\right)=0 \Longrightarrow \triangle_{B_1,A_1}^{m_1}\left(\delta_{B_2,A_2}^n(X)\right)=0 \ {\rm for \ all}\ m_1 \geq m,
$$
if either of the hypotheses (c), or, $(a)\wedge (d)$, or, $(b)\wedge (d)$ is satisfied, then
$$
\left((SB_1,TA_1), (B_2,A_2)\right)\in left-\left(X,(m+t-1,n)\right)-symmetric.
$$
\end{demo}

A particularly interesting case of Proposition \ref{pro1} is obtained upon choosing
$$
A_1=A_2=A, B_1=B_2=A^*, S_1=S_2=B^* \ {\rm and}\ T_1=T_2=B.
$$

We have:

\begin{cor}\label{cor01} If $[A,B]=[A,B^*]=0,\ A\in (m,n)$-isosymmetric and $B$ is both left-$r$-invertible and $s$-symmetric, then
$$
AB\in left-(m+r-1,n+s-1)-symmetric.
$$
\end{cor}

A generalization of Proposition \ref{pro1} is obtained upon replacing the pair of operators $(S,T)$ by the pairs $(S_i,T_i),\ i=1,2$. Observe that the hypothesis
$$
(S_1,T_1)\in left-(X,t)-invertible \Longrightarrow \delta_{E,F}^j\left(\triangle_{S_1,T_1}^t(X)\right)=\triangle_{E,F}^k\left(\triangle_{S_1,T_1}^t(X)\right)=0
$$
for all $1\leq j, k$ and $(E,F)=(S_2,T_2)\vee (B_i,A_i),\ i=1,2,$ and the hypothesis
$$
(S_2,T_2)\in left-(X,s)-symmetric \Longrightarrow \delta_{E,F}^j\left(\delta_{S_2,T_2}^s(X)\right)=\triangle_{E,F}^k\left(\delta_{S_2,T_2}^s(X)\right)=0
$$
for all $1\leq j, k$ and $(E,F)=(S_1,T_1)\vee (B_i,A_i),\ i=1,2.$

\begin{cor}\label{cor02} If $[A_i,T_j]=[B_i,S_j]=[S_1,S_2]=[B_1,B_2]=[A_1,A_2]=0, \ 1 \leq i,j\leq 2$,

\vskip4pt\noindent (i) $\left((B_1,A_1), (B_2,A_2)\right)\in$ left-$\left(X,(m,n)\right)$-symmetric;

\vskip4pt\noindent (ii) $(S_1,T_1)\in$ left-$\left(X,r)\right)$-invertible; and

\vskip4pt\noindent (iii) $(S_2,T_2)\in(X,s)$-symmetric,

\noindent then
$$
\left((B_1S_1,A_1T_1), (B_2S_2,A_2T_2)\right)\in left-\left(X,(m+r-1,n+s-1)\right)-symmetric.
$$
\end{cor}

\begin{demo} In view of the commutativity hypotheses, Proposition \ref{pro1} implies that if (i) and (ii) are satisfied, then

\vskip4pt\noindent (iv) $\left((B_1S_1,A_1T_1), (B_2,A_2)\right)\in$ left-$\left(X,(m+r-1,n)\right)$-symmetric;

\noindent another application of Proposition \ref{pro1}, this time since (iv) and (iii) are satisfied, now implies the result.
\end{demo}

Theorem \ref{thm1}, which we now prove, is a generalization of Corollary \ref{cor02}.

\vskip8pt\noindent{\bf Proof of Theorem \ref{thm1}.} Defining the positive integers $m,n,r$ and $s$ as in the statement of the theorem, the commutativity hypotheses imply
\begin{eqnarray*} \delta_{S_2B_2,T_2A_2}^{n+s-1}&=&\left\{L_{S_2}\left(L_{B_2}-R_{A_2}\right)+\left(L_{S_2}-R_{T_2}\right)R_{A_2}\right\}^{n+s-1}\\
&=&  \sum_{k=0}^{n+s-1}\left(\begin{array}{clcr}n+s-1\\k\end{array}\right) L_{S_2}^{n+s-1-k}R_{A_2}^k\delta_{B_2,A_2}^{n+s-1-k}\delta_{S_2,T_2}^k;
\end{eqnarray*}

\begin{eqnarray*} \triangle_{S_1B_1,T_1A_1}^{m+r-1}&=&\left\{L_{S_1}\left(L_{B_1}L_{A_1}-I\right)R_{T_1}+\left(L_{S_1}R_{T_1}-I\right)\right\}^{m+r-1}\\
&=&  \sum_{j=0}^{m+r-1}\left(\begin{array}{clcr}m+r-1\\j\end{array}\right) \left(L_{S_1}R_{T_1}\right)^{m+r-1-j}\triangle_{B_1,A_1}^{m+r-1-j}\triangle_{S_1,T_1}^j;
\end{eqnarray*}
and
$$
\triangle_{S_1B_1,T_1A_1}^{m+r-1}\left(\delta_{S_2B_2,T_2A_2}^{n+s-1}(X)\right)=\left\{\sum_{j=0}^{m+r-1}\sum_{k=0}^{n+s-1}\left(\begin{array}{clcr}m+r-1\\j\end{array}\right)
\left(\begin{array}{clcr}n+s-1\\k\end{array}\right)\times\right.
$$
$$
\left.\times\left(L_{S_1}R_{T_1}\right)^{m+r-1-j}L_{S_2}^{n+s-1-k}R_{A_2}^k \triangle_{B_1,A_1}^{m+r-1-j}\triangle_{S_1,T_1}^j\delta_{B_2,A_2}^{n+s-1-k}\delta_{S_2,T_2}^k\right\}(X)
$$
where the commutativity hypotheses guarantee the commutativity of all the operator entries within the curly brackets. To complete the proof we observe now that:

if $j\geq r$ and $k\geq s$, then $\triangle_{S_1,T_1}^j\left(\delta_{S_2,T_2}^k(X)\right)=0$;

if $j\leq r-1$ and $k\geq s$, then $\triangle_{B_1,A_1}^{m+r-1-j}\left(\delta_{S_2,T_2}^k(X)\right)=0$;

if $j\geq r$ and $k\leq s-1$, then $\triangle_{S_1,T_1}^j\left(\delta_{B_2,A_2}^{n+s-1-k}(X)\right)=0$

and (finally)

if $j\leq r-1$ and $k\leq s-1$, then $\triangle_{B_1,A_1}^{m+r-1-j}\left(\delta_{B_2,A_2}^{n+s-1-k}(X)\right)=0$.

\

\section {\sfstp Tensor Products.} Let $\H\otimes\H$, endowed with a reasonable uniform cross norm, denote the completion of the algebraic tensor product of $\H$ with itself, and let $A\otimes B$ denote the tensor product operator defined by $A,B\in \B$. Theorem \ref{thm1} has applications to tensor products.

Consider operators $A_i,B_i, \ i=1,2$. Since
\begin{eqnarray*}
& &\triangle_{B_1\otimes I,A_1\otimes I}^m\left(\delta_{B_2\otimes I,A_2\otimes I}^{n}\right)\\
&=&\sum_{j=0}^{m}\sum_{k=0}^{n}(-1)^{j+k}\left(\begin{array}{clcr}m\\j\end{array}\right)
\left(\begin{array}{clcr}n\\k\end{array}\right)
\left(L_{B_1\otimes I}R_{A_1\otimes I}\right)^{m-j}\left(L_{B_2\otimes I}^{n-k}R_{A_2\otimes I}^k\right)\\
&=&\sum_{j=0}^{m}\sum_{k=0}^{n}(-1)^{j+k}\left(\begin{array}{clcr}m\\j\end{array}\right)
\left(\begin{array}{clcr}n\\k\end{array}\right)
\left(L_{{B_1}^{m-j}\otimes I}R_{{A_1}^{m-j}\otimes I}\right)\left(L_{{B_2}^{n-k}\otimes I}R_{{A_2}^k\otimes I}\right),
\end{eqnarray*}
$$
\triangle_{B_1\otimes I,A_1\otimes I}^m\left(\delta_{B_2\otimes I,A_2\otimes I}^{n}(X\otimes X)\right)=\triangle_{B_1,A_1}^m\left(\delta_{B_2,A_2}^{n}(X)\right)\otimes X.
$$
Hence, if $\left((B_1,A_1), (B_2,A_2)\right)\in$ left-$\left(X,(m,n)\right)$-symmetric, then
$$
\triangle_{B_1\otimes I,A_1\otimes I}^m\left(\delta_{B_2\otimes I,A_2\otimes I}^{n}(X\otimes X)\right)=0
$$
(and, arguing similarly,
$$
\left. \triangle_{I\otimes B_1,I\otimes A_1}^m\left(\delta_{I\otimes B_2,I\otimes A_2}^{n}(X\otimes X)\right)=0\right).
$$

\begin{cor}\label{cor03} Let $E_i, F_i, P_i, Q_i \in \B,\ i=1,2$, be such that $[E_1,E_2]=[F_1,F_2]=0,$
$$
\triangle_{E_1,F_1}^{m_1}\left(\delta_{E_2,F_2}^{n_1}(X)\right)=0=\triangle_{P_1,Q_1}^{r_1}\left(\delta_{E_2,F_2}^{n_2}(X)\right), \ {\rm and}
$$
$$
\triangle_{E_1,F_1}^{m_2}\left(\delta_{P_2,Q_2}^{s_1}(X)\right)=0=\triangle_{P_1,Q_1}^{r_2}\left(\delta_{P_2,Q_2}^{s_2}(X)\right).
$$
Then
$$
\left((E_1\otimes P_1, F_1\otimes Q_1), (E_2\otimes P_2, F_2\otimes Q_2)\right)\in left-\left(X\otimes X, (m+r-1,n+s-1)\right)-symmetric,
$$
where $m=\max (m_1,m_2), n=\max(n_1,n_2), r=\max(r_1,r_2)$ and $s=\max(s_1,s_2)$.
\end{cor}

\begin{demo} If we let
$$
B_i=E_i\otimes I, A_i=F_i\otimes I, S_i=I\otimes P_i \ {\rm  and}   \ T_i=I\otimes Q_i,\ i=1,2,
$$
then the hypotheses of Theorem \ref{thm1} are satisfied (with $X\otimes X$ playing the role of $X$). Hence, since
$$
S_iB_i=E_i\otimes P_i \ {\rm and} \ T_iA_i=F_i\otimes Q_i, \ i=1,2,
$$
the proof follows.
\end{demo}

For commuting $n$-isometries $S,T\in \B$, Corollary \ref{cor03} take the form:

\begin{cor}\label{cor04} If $[S,T]=0$ and
$$
\triangle_{S^*,S}^n\left(\delta_{S^*,S}^n(I)\right)=\triangle_{T^*,T}^n\left(\delta_{T^*,T}^n(I)\right)
=\triangle_{T^*,T}^n\left(\delta_{S^*,S}^n(I)\right)=\triangle_{S^*,S}^n\left(\delta_{T^*,T}^n(I)\right)=0,
$$
then
$$
\triangle_{S^*\otimes T^*,S\otimes T}^{2n-1}\left(\delta_{S^*\otimes T^*,S\otimes T}^{2n-1}(I\otimes I)\right)=0.
$$
\end{cor}

\

\section {\sfstp Perturbation by commuting nilpotents.} It is well known, see for example \cite{TL}, that if $d_{B,A}=\delta_{B,A}\vee \triangle_{B,A}$, $d_{B,A}^m(I)=0$ and $N$ is an $n$-nilpotent which commutes with $A$, then $d_{B,A+N}^{m+n-1}(I)=0$. This extends to $A,B$ such that $d_{B,A}^m(X)=0$ for some $X\in \B$, as the following argument shows:
\begin{eqnarray*}
\triangle_{B,A+N}^{m+n-1}(X)&=&\left\{\left(L_{B}R_{A}-I\right)+L_B R_N\right\}^{m+n-1}(X)\\
&=& \left\{ \sum_{j=0}^{m+n-1}\left(\begin{array}{clcr}m+n-1\\j\end{array}\right)\left(L_{B}R_N\right)^{j}\triangle_{B,A}^{m+n-1-j}\right\}(X)\\
&=& 0,
\end{eqnarray*}
since $R_N^j=0$ for all $j\geq n$ and $\triangle_{B,A}^{m+n-1-j}(X)=0$ for all $j\leq n-1 (\Longrightarrow m+n-1-j\geq m)$; again
\begin{eqnarray*}
\delta_{B,A+N}^{m+n-1}(X)&=&\left\{\left(L_{B}-R_{A}\right)- R_N\right\}^{m+n-1}(X)\\
&=& \left\{ \sum_{j=0}^{m+n-1}(-1)^j\left(\begin{array}{clcr}m+n-1\\j\end{array}\right) R^j_N \delta_{B,A}^{m+n-1-j}\right\}(X)\\
&=& 0.
\end{eqnarray*}
for $R_N^j=0$ for all $j\geq n$ and $\delta_{B,A}^{m+n-1-j}(X)=0$ for all $j\leq n-1$.

\

This argument extends to perturbation by commuting nilpotents of operators $\left((B_1,A_1),(B_2,A_2)\right)\in {\rm left}-(X, (m,n))-{\rm symmetric}$.

\

\noindent{\bf Proof of Theorem \ref{thm2}.} The commutativity hypotheses $[A_1,A_2]=[B_1,B_2]=0$ implies
$$
\delta_{B_2,A_2}^n\left(\triangle_{B_1,A_1}^m(X)\right)=\triangle_{B_1,A_1}^m\left(\delta_{B_2,A_2}^n(X)\right)=0.
$$
Set
$$
\delta_{B_2,A_2}^n(X)=Y.
$$
Then $\triangle_{B_1,A_1}^m(Y)=0$, and it follows from the argument above that if $[A_1,M_1]=0$ and $M_1^{m_1}=0$, then
$$
\triangle_{B_1,A_1+M_1}^{m+m_1-1}(Y)=0 \Longleftrightarrow \triangle_{A_1^*+M_1^*,B_1^*}^{m+m_1-1}(Y^*)=0,
$$
and hence if $[B_1,N_1]=0, N_1^{n_1}=0$, then
$$
\triangle_{A_1^*+M_1^*,B_1^*+N_1^*}^{m+m_1+n_1-2}(Y^*)=0 \Longleftrightarrow \triangle_{B_1+N_1, A_1+M_1}^{m+m_1+n_1-2}(Y)=0.
$$
Let now $\triangle_{B_1+N_1, A_1+M_1}^{m+m_1+n_1-1}(X)=Z$. Then, upon assuming the full force of the commutativity hypotheses, we have
\begin{eqnarray*}
\delta_{B_2,A_2}^n(Z)=0 &\Longrightarrow& \delta_{B_2,A_2+M_2}^{n+m_2-1}(Z)=0 \Longleftrightarrow \delta_{A_2^*+M_2^*,B_2^*}^{n+m_2-1}(Z^*)=0\\
 &\Longrightarrow& \delta_{A_2^*+M_2^*,B_2^*+N_2^*}^{n+m_2+n_2-2}(Z^*)=0 \Longleftrightarrow \delta_{B_2+N_2,A_2+M_2}^{n+m_2+n_2-2}(Z)=0\\
 &\Longleftrightarrow& \delta_{B_2+N_2,A_2+M_2}^{n+m_2+n_2-2}\left(\triangle_{B_1+N_1,A_1+M_1}^{n+m_2+n_2-2}(X)\right)\\
  & & = \triangle_{B_1+N_1,A_1+M_1}^{n+m_2+n_2-2}\left(\delta_{B_2+N_2,A_2+M_2}^{n+m_2+n_2-2}(X)\right)=0. 
\end{eqnarray*} 
We leave it to the reader to verify that if one drops the hypotheses $[A_1,A_2]=[B_1,B_2]=[M_1,M_2]=[N_1,N_2]=0$ in the statement of Theorem \ref{thm2}, then
$$
\delta_{B_2,A_2}^n\left(\triangle_{B_1,A_1}^m(X)\right)=0 \Longrightarrow \delta_{B_2+N_2,A_2+M_2}^{n+m_2+n_2-2}\left(\triangle_{B_1+N_1,A_1+M_1}^{n+m_2+n_2-2}(X)\right)=0
$$
only. For $(X,(m,n))$-isosymmetric operators $A$, Theorem \ref{thm2} implies

\begin{cor}\label{cor05} If $[A,N]=0=N^{n_1}$ and $A\in (X,(m,n))$-isosymmetric, then $A+N \in (X,(m+2n_1-2,n+2n_1-2))$-isosymmetric.
\end{cor}

\

\section {\sfstp Application to Drazin invertible operators.} Let $T\in \B$ be a Drazin invertible operator. Then there exists a decomposition
$$
\H=\H_1 \oplus \H_2
$$
of $\H$ and a decomposition
$$
T=T_1\oplus T_2 \ (=T\mid_{\H_1}\oplus T\mid_{\H_2})
$$
of $T$ such that $T_1$ is invertible and $T_2^p=0$ for some integer $p\geq 1$. The Drazin inverse $T_d$ of $T$ has a representation
$$
T_d=T_1^{-1}\oplus 0\in B(\H_1 \oplus \H_2)
$$
\cite[Theorem 2.23]{DR}. In the following we use this representation of Drazin invertible operators to prove Theorem \ref{thm3}.

\

\noindent {\bf Proof of Theorem \ref{thm3}}  Let $X\in B(\H_1\oplus \H_2)$ have the matrix representation $\left[X_{i,j}\right]_{i,j=1}^2$. Then

\begin{eqnarray*}
0&=&\triangle_{T^*,T}^m\left(\delta_{T^*,T}^n(X)\right)=\delta_{T^*,T}^n\left(\triangle_{T^*,T}^m(X)\right)\\
\noindent (1) \hspace{2cm}&=&\sum_{l=0}^n\sum_{k=0}^m(-1)^{l+k}\left(\begin{array}{clcr}n\\l\end{array}\right)
\left(\begin{array}{clcr}m\\k\end{array}\right)\left\{{T^*}^{(n-l+m-k)}X T^{m-k+l}\right\}\\
&=&\sum_{l=0}^n\sum_{k=0}^m(-1)^{l+k}\left(\begin{array}{clcr}n\\l\end{array}\right)
\left(\begin{array}{clcr}m\\k\end{array}\right)\left[{T_i^*}^{(n-l+m-k)}X_{ij} T_j^{m-k+l}\right]_{i,j=1}^2;
\end{eqnarray*}

\begin{eqnarray*}
0&=&\triangle_{T_d^*,T}^m\left(\delta_{T^*,T}^n(X)\right)=\delta_{T^*,T}^n\left(\triangle_{T_d^*,T}^m(X)\right)\\
\noindent (2) \hspace{0.8cm}&=&\sum_{l=0}^n\sum_{k=0}^m(-1)^{l+k}\left(\begin{array}{clcr}n\\l\end{array}\right)
\left(\begin{array}{clcr}m\\k\end{array}\right)\left(\begin{array}{cccc}{T_1^*}^{(n-l-m+k)}X_{11} T_1^{m-k+l}&{T_1^*}^{(n-l-m+k)}X_{12} T_2^{m-k+l}\\0&0\end{array}\right)
\end{eqnarray*}
and
\begin{eqnarray*}
0&=&\triangle_{T^*,T}^m\left(\delta_{T_d^*,T}^n(X)\right)=\delta_{T_d^*,T}^n\left(\triangle_{T^*,T}^m(X)\right)\\
\noindent (3) \hspace{0.8cm}&=&\sum_{l=0}^n\sum_{k=0}^m(-1)^{l+k}\left(\begin{array}{clcr}n\\l\end{array}\right)
\left(\begin{array}{clcr}m\\k\end{array}\right)\left(\begin{array}{cccc}{T_1^*}^{(m-k-n+l)}X_{11} T_1^{m-k+l}&{T_1^*}^{(m-k-n+l)}X_{12} T_2^{m-k+l}\\0&0\end{array}\right).
\end{eqnarray*}
We prove that $X_{12}=0$. Letting $s$ denote either of $n-l$ and $-n+l$, we have
\begin{eqnarray*}
& &\sum_{l=0}^n\sum_{k=0}^m(-1)^{l+k}\left(\begin{array}{clcr}n\\l\end{array}\right)
\left(\begin{array}{clcr}m\\k\end{array}\right){T_i^*}^{-m+k+s}X_{12} T_2^{m-k+l}=0\\
&\Longleftrightarrow& \sum_{l=0}^n(-1)^{l}\left(\begin{array}{clcr}n\\l\end{array}\right){T_1^*}^s\left\{\sum_{k=0}^{m-1}(-1)^{k}\left(\begin{array}{clcr}m\\k\end{array}\right){T_1^*}^{-m+k}X_{12} T_2^{m-k}+X_{12}\right\}T_2^l=0\\
&\Longrightarrow& \left\{\left(\sum_{l=1}^n(-1)^{l}\left(\begin{array}{clcr}n\\l\end{array}\right)\sum_{k=0}^{m-1}(-1)^{k}\left(\begin{array}{clcr}m\\k\end{array}\right){T_1^*}^{s-m+k}X_{12} T_2^{m-k+l}+{T_1^*}^sX_{12}T_2^l\right)\right.\\
& & \left.+\sum_{k=0}^{m-1}(-1)^{k}\left(\begin{array}{clcr}m\\k\end{array}\right){T_1^*}^{s-m+k}X_{12} T_2^{m-k}+{T_1^*}^sX_{12}\right\}T_2^{p-1}=0\\
&\Longleftrightarrow& {T_1^*}^sX_{12}T_2^{p-1}=0 \Longleftrightarrow X_{12}T_2^{p-1}=0.
\end{eqnarray*}
Repeating this argument next by multiplying the expression within the curly brackets by $T_2^{p-2}$ on the right, and then so on, it is seen that
$$
X_{12}T_2^r=0 \ {\rm for \ all}\ r=1,2,\cdots.
$$
Consequently, it follows (from the expression within the curly brackets) that $X_{12}=0$. A similar argument shows that $X_{21}=0$ (in (1)).

\

\noindent (i) Equality (1) implies $T_1 \in (X_{11},(m,n))$-isosymmetric, and this (since $T_1$ is invertible) implies by Lemma \ref{lem0} that $T_1^{-1}\in (m,n)$-isosymmetric, i.e.,
\begin{eqnarray*}
& &\sum_{l=0}^n\sum_{k=0}^m(-1)^{l+k}\left(\begin{array}{clcr}n\\l\end{array}\right)
\left(\begin{array}{clcr}m\\k\end{array}\right)\left\{{T_1^*}^{-(n-l+m-k)}X_{11} T_1^{-m+k-l}\right\}=0\\
&\Longleftrightarrow& \sum_{l=0}^n(-1)^{l}\left(\begin{array}{clcr}n\\l\end{array}\right)
{T_1^*}^{-(n-l)}\left\{\sum_{k=0}^m(-1)^{k}\left(\begin{array}{clcr}m\\k\end{array}\right){T_1^*}^{-(m-k)}X_{11} T_1^{-(m-k)}\right\}T_1^{-l}=0\\
&\Longleftrightarrow& \sum_{l=0}^n(-1)^{l}\left(\begin{array}{clcr}n\\l\end{array}\right)
{T_1^*}^{-(n-l)}\left\{\sum_{k=0}^m(-1)^{k}\left(\begin{array}{clcr}m\\k\end{array}\right){T_1^*}^{-(m-k)}X_{11} T_1^{-(m-k)}\right\}T_1^{n-l}=0\\
&\Longleftrightarrow& \triangle_{{T_1^*}^{-1},T_1}^n\left(\triangle_{{T_1^*}^{-1},T_1^{-1}}^m(X_{11})\right)=0\\
&\Longleftrightarrow& \triangle_{{T_d^*},T_1}^n\left(\triangle_{{T_d^*},T_d}^m(X)\right)=0.
\end{eqnarray*}
Again, $T_1^{-1}\in (X_{11},(m,n))$-isosymmetric if and only if
\begin{eqnarray*}
& &\sum_{k=0}^m(-1)^{k}\left(\begin{array}{clcr}m\\k\end{array}\right) T_1^{*(-m+k)}
\left\{\sum_{l=0}^n(-1)^{l}\left(\begin{array}{clcr}n\\l\end{array}\right){T_1^*}^{-(n-l)}X_{11} T_1^{-l}\right\}T_1^{-(m-k)}=0\\
&\Longleftrightarrow& \sum_{k=0}^m(-1)^{k}\left(\begin{array}{clcr}m\\k\end{array}\right) T_1^{*(-m+k)}
\left\{\sum_{l=0}^n(-1)^{l}\left(\begin{array}{clcr}n\\l\end{array}\right){T_1^*}^{-(n-l)}X_{11} T_1^{-l}\right\}T_1^k=0\\
&\Longleftrightarrow& \delta_{{T_1^*}^{-1},T_1}^m\left(\delta_{{T_1^*}^{-1},T_1^{-1}}^n(X_{11})\right)=0\\
&\Longleftrightarrow& \delta_{{T_d^*},T_1}^m\left(\delta_{{T_d^*},T_d}^n(X)\right)=0.
\end{eqnarray*}

\

\noindent (ii) Since $X_{12}=0$, and
\begin{eqnarray*}
& &\sum_{l=0}^n\sum_{k=0}^m(-1)^{l+k}\left(\begin{array}{clcr}n\\l\end{array}\right)
\left(\begin{array}{clcr}m\\k\end{array}\right)\left\{{T_1^*}^{(-m+k+n-l)}X_{11} T_1^{m-k+l}\right\}=0\\
&\Longleftrightarrow& \sum_{l=0}^n(-1)^{l}\left(\begin{array}{clcr}n\\l\end{array}\right)
{T_1^*}^{-(n-l)}\left\{\sum_{k=0}^m(-1)^{k}\left(\begin{array}{clcr}m\\k\end{array}\right){T_1^*}^{m-k}X_{11} T_1^{-(m-k)}\right\}T_1^{-l}=0\\
&\Longleftrightarrow& \sum_{l=0}^n(-1)^{l}\left(\begin{array}{clcr}n\\l\end{array}\right)
{T_1^*}^{-(n-l)}\left\{\sum_{k=0}^m(-1)^{k}\left(\begin{array}{clcr}m\\k\end{array}\right){T_1^*}^{m-k}X_{11} T_1^{-(m-k)}\right\}T_1^{n-l}=0\\
&\Longleftrightarrow& \triangle_{{T_1^*}^{-1},T_1}^n\left(\triangle_{{T_1^*},T_1^{-1}}^m(X)\right)=0,\ {\rm we\ have\ that}
\end{eqnarray*}
$\triangle_{{T_d^*},T}^n\left(\triangle_{{T^*},T_d}^m(X)\right)=0.$

\

\noindent (iii) In this case
\begin{eqnarray*}
& &\sum_{l=0}^n\sum_{k=0}^m(-1)^{l+k}\left(\begin{array}{clcr}n\\l\end{array}\right)
\left(\begin{array}{clcr}m\\k\end{array}\right)\left\{{T_1^*}^{(m-k-n+l)}X_{11} T_1^{m-k+l}\right\}=0\\
&\Longleftrightarrow& \sum_{l=0}^n\sum_{k=0}^m(-1)^{l+k}\left(\begin{array}{clcr}n\\l\end{array}\right)
\left(\begin{array}{clcr}m\\k\end{array}\right)\left\{{T_1^*}^{-(m-k-n+l)}X_{11} T_1^{-(m-k+l)}\right\}=0\\
&\Longleftrightarrow& \sum_{l=0}^n(-1)^{l}\left(\begin{array}{clcr}n\\l\end{array}\right)
{T_1^*}^{(n-l)}\left\{\sum_{k=0}^m(-1)^{k}\left(\begin{array}{clcr}m\\k\end{array}\right){T_1^*}^{-(m-k)}X_{11} T_1^{-(m-k)}\right\}T_1^{-l}=0\\
&\Longleftrightarrow& \sum_{l=0}^n(-1)^{l}\left(\begin{array}{clcr}n\\l\end{array}\right)
{T_1^*}^{(n-l)}\left\{\sum_{k=0}^m(-1)^{k}\left(\begin{array}{clcr}m\\k\end{array}\right){T_1^*}^{-(m-k)}X_{11} T_1^k\right\}T_1^{-l}=0\\
&\Longleftrightarrow& \delta_{{T_1^*},T_1^{-1}}^n\left(\delta_{{T_1^*}^{-1},T_1}^m(X_{11})\right)=0\\
&\Longleftrightarrow& \delta_{{T^*},T_d}^n\left(\delta_{{T_d^*},T}^m(X)\right)=0.
\end{eqnarray*}


\vskip10pt \noindent\normalsize\rm B.P. Duggal,{University of Ni\v s,
Faculty of Sciences and Mathematics,
P.O. Box 224, 18000 Ni\v s, Serbia}.

\noindent\normalsize \tt e-mail:  bpduggal@yahoo.co.uk

\vskip6pt\noindent \noindent\normalsize\rm I. H. Kim, Department of
Mathematics, Incheon National University, Incheon, 22012, Korea.\\
\noindent\normalsize \tt e-mail: ihkim@inu.ac.kr

\end{document}